\documentclass[12pt]{amsart}

\usepackage{amsmath,amsthm}
\usepackage{amssymb,latexsym}
\usepackage{graphicx}
\theoremstyle{definition}

\theoremstyle{remark}

\numberwithin{equation}{section}

\newcommand{\bee}{\begin{equation}}
\newcommand{\ene}{\end{equation}}

\newcommand{\al}{\alpha}

\begin{document}

\title[A Diffusion Problem with Neumann Boundary Control]{A Diffusion Problem with Neumann Boundary Control Utilizing Total Mass}

\author{M. Salman}
\address{Department of Mathematical Science, University of Delaware, Newark, DE 19716} \email{msalmanz@gmail.com}



\date{\today}

\begin{abstract}
The author studies the diffusion problem $u_t=u_{xx},\ 0<x<1,\ t>0;
\ u(x,0)=0,$ and $-u_x(0,t)=u_x(1,t)=\phi(t),$ where $\phi(t)$ is a
control function that ensures that the total mass  $\int_0^1
u(x,t_k)dx$ stays between two predetermined values. Mathematically,
$\phi(t)=1$ for $t_{2k} < t<t_{2k+1}$ and $\phi(t)=-1$ for $t_{2k+1}
<t<t_{2k+2},\ k=0,1,2,\ldots$ with $t_0=0$ and the sequence $t_{k}$
is determined by the equations $\int_0^1 u(x,t_k)dx = M,$ for
$k=1,3,5,\dots,$ and $\int_0^1 u(x,t_k)dx = m,$ for $k=2,4,6,\dots$
and where $0<m<M<u_0$. Note that the switching times $t_k$ are
unknowns. Determination of switching times $t_k$ and their
differences $t_{k+1}-t_k$ are obtained.  Numerical verifying
examples are presented.

\end{abstract}
\maketitle


\numberwithin{equation}{section}








\section{Introduction}

 As motivation for the mathematical problems considered
in this work, consider a chamber in the form of a long linear
transparent tube. We allow for the introduction or removal of
material in a gaseous state at the ends of the tube. The material
diffuses throughout the tube with or without reaction with other
materials. By illuminating the tube on one side with a light source
with a frequency range spanning the absorption range for the
material and collecting the residual light that passes through the
tube with photo-reception equipment, we can obtain a measurement of
the total mass of material contained in the tube as a function of
time. Using the total mass as switch points for changing the
boundary conditions for introduction or removal of material. The
objective is to keep the total mass of material in the tube
oscillating between two set values such as $m<M$. The physical
application for such a system is the control of reaction diffusion
systems such as production of a chemical material in a reaction
chamber via the introduction of reactants at the boundary of
chamber.

\section{Definition of the Problem}

We consider the diffusion problem
\begin{equation}
\begin{array}{llll}
u_{t} & = & u_{xx}, & \,\,\,\,\,\,\,\,\,\,\,\,\,\,\,\,\,\,\,\,\,\,\,\,\,\,\,%
\,\,\,\,\,\,\,\,t>0,\,\,\,0<x<1, \\
u(x,0) & = & 0 & \,\,\,\,\,\,\,\,\,\,\,\,\,\,\,\,\,\,\,\,\,\,\,\,\,\,\,\,\,%
\,\,\,\,\,0<x<1, \\
u_{x}(0,t) & = & -\phi (t) & \,\,\,\,\,\,\,\,\,\,\,\,\,\,\,\,\,\,\,\,\,\,\,%
\,\,\,\,\,\,\,\,\,\,\,\,t>0, \\
u_{x}(1,t) & = & \phi (t) & \,\,\,\,\,\,\,\,\,\,\,\,\,\,\,\,\,\,\,\,\,\,\,\,%
\,\,\,\,\,\,\,\,\,\,\,t>0.
\end{array}
\label{5.1}
\end{equation}
where the flux $\phi (t)$ is set in such a way to control the level
of the total mass $\mu (t)=\int_{0}^{1}u(x,t)dx$. Namely,
\begin{equation*}
\phi (t)=\left\{
\begin{array}{ll}
1 & \,\,\,\,\,\,\,\,\,\,\text{if}\,\,\,\,\,t_{2n}<t<t_{2n+1} \\
&  \\
-1 & \,\,\,\,\,\,\,\,\,\,\text{if}\,\,\,\,\,t_{2n+1}<t<t_{2n+2}
\end{array}
\right.
\end{equation*}
where the sequence
\begin{equation*}
0=t_{0}<t_{1}<t_{2}<\ldots ,
\end{equation*}
satisfies
\begin{equation*}
\begin{array}{lll}
\mu (t_{2n}) & = & m,\,\,\,\,\,\,\,\,\,\,\,\,\,\,\,\,\,\,\,\,\,\,\,\,%
\,n=1,2,3,\ldots \\
\mu (t_{2n+1}) & = & M,\,\,\,\,\,\,\,\,\,\,\,\,\,\,\,\,\,\,\,\,\,\,\,\,%
\,n=0,1,2,\ldots
\end{array}
\end{equation*}
with $m$ and $M$ are certain positive threshold to ensure $m<\mu
<M.$ The total mass has a time derivative equal to
\begin{eqnarray}
\mu ^{\prime }(t) &=&\int_{0}^{1}u_{t}(x,t)dx  \label{5.2} \\
&=&\int_{0}^{1}u_{xx}(x,t)dx  \notag \\
&=&u_{x}(1,t)-u_{x}(0,t)  \notag \\
&=&2\phi (t).  \notag
\end{eqnarray}
where we use the fact that $u(x,t)$ satisfies (\ref{5.1}). This relationship
will allow us to explicitly find the time switches $t_{n}.$

Since the flux $\phi (t)$ is initially set at $1,$ and $\mu (0)=0$ due to
the initial condition on $u,$ then
\begin{equation*}
\mu (t)=2t,\,\,\,\,\,\,\,\,\,\,\,\,\,\,\,\,\,\,\,\,\,\,\,\,\,\,\,\,\,\,0\le
t\le t_{1}.
\end{equation*}
where $t_{1}$ can be calculated as a solution of $\mu (t_{1})=M,$ which
implies $t_{1}=M/2.$

At the second time interval $t_{1}<t<t_{2},$ the flux is reversed, that is $%
\phi (t)=-1,$ therefore the total mass will be
\begin{equation*}
\mu
(t)=-2t+M+2t_{1},\,\,\,\,\,\,\,\,\,\,\,\,\,\,\,\,\,\,\,\,\,\,\,\,\,\,\,\,\,%
\,t_{1}\le t\le t_{2},
\end{equation*}
where due to the continuity of $\mu (t)$ at\thinspace $t_{1}.$ The second
time switch $t_{2}$ will be calculated as a solution of $\mu (t_{2})=m,$
this gives
\begin{eqnarray*}
t_{2} &=&t_{1}+(\frac{M}{2}-\frac{m}{2}) \\
&=&M-\frac{m}{2},
\end{eqnarray*}
Now, if we bring the flux back to $\phi (t)=1,$ where $\,t_{2}\le t\le
t_{3}, $ we will find
\begin{equation*}
\mu
(t)=2t+m-2t_{2},\,\,\,\,\,\,\,\,\,\,\,\,\,\,\,\,\,\,\,\,\,\,\,\,\,\,\,\,\,%
\,t_{2}\le t\le t_{3},
\end{equation*}
Setting $\mu (t_{3})=M$ yields
\begin{eqnarray*}
t_{3} &=&t_{2}+(\frac{M}{2}-\frac{m}{2}) \\
&=&\frac{3M}{2}-m,
\end{eqnarray*}
Following the same pattern, we will successively get
\begin{eqnarray*}
\mu (t)
&=&-2t+M+2t_{3},\,\,\,\,\,\,\,\,\,\,\,\,\,\,\,\,\,\,\,\,\,\,\,\,\,\,\,\,\,%
\,t_{3}\le t\le t_{4}, \\
&& \\
t_{4} &=&t_{3}+(\frac{M}{2}-\frac{m}{2}) \\
&=&2M-\frac{3m}{2}, \\
&& \\
\mu(t)&=&2t+m-2t_{4},\,\,\,\,\,\,\,\,\,\,\,\,\,\,\,\,\,\,\,\,\,\,\,\,\,\,\,%
\,\, \,t_{4}\le t\le t_{5}, \\
&& \\
t_{5} &=&t_{4}+(\frac{M}{2}-\frac{m}{2}) \\
&=&\frac{5M}{2}-2m, \\
&\vdots & \\
\mu (t)&=&\left\{
\begin{array}{ll}
2t+m-2t_{2n} & \,\,\,\,\,\,\,\,\,\,\text{if}\,\,\,\,\,t_{2n}<t<t_{2n+1} \\
&  \\
-2t+M+2t_{2n+1} & \,\,\,\,\,\,\,\,\,\,\text{if}\,\,\,\,\,t_{2n+1}<t<t_{2n+2}
\end{array}
\right. \\
&& \\
t_{n} &=&t_{n-1}+(\frac{M}{2}-\frac{m}{2}) \\
&=&\frac{nM}{2}-\frac{(n-1)m}{2},
\end{eqnarray*}
where $n\ge 1.$ We infer from these calculations that the time switches $%
t_{n}$ are equispaced with $t_{n}-t_{n-1}=0.5M-0.5m.$

\section{Finite difference analysis}

We study a finite difference discretization for the control problem (\ref
{5.1}). We consider the implicit backward scheme
\begin{eqnarray*}
\frac{U_{j}^{n+1}-U_{j}^{n}}{\Delta t} &=&\frac{%
U_{j-1}^{n+1}-2U_{j}^{n+1}+U_{j+1}^{n+1}}{(\Delta x)^{2}},\,\,\,\,\,\,\,\,\,%
\,j=1,\dots ,J-1,\,\,\,n=0,1,2,\dots \\
U_{j}^{0}
&=&0,\,\,\,\,\,\,\,\,\,\,\,\,\,\,\,\,\,\,\,\,\,\,\,\,\,\,\,\,\,\,\,\,\,\,\,%
\,\,\,\,\,\,\,\,\,\,\,\,\,\,\,\,j=0,\dots ,J \\
\frac{U_{1}^{n+1}-U_{0}^{n+1}}{\Delta x} &=&-\phi (\tau
_{n+1}),\,\,\,\,\,\,\,\,\,\,\,\,\,\,\,\,\,\,\,\,\,\,\,\,\,\,\,\,n=0,1,2,\dots
\\
\frac{U_{J}^{n+1}-U_{J-1}^{n+1}}{\Delta x} &=&\phi (\tau
_{n+1}),\,\,\,\,\,\,\,\,\,\,\,\,\,\,\,\,\,\,\,\,\,\,\,\,\,\,\,\,\,\,\,\,%
\,n=0,1,2,\dots
\end{eqnarray*}
where $\tau _{n}$ is not necessarily a uniform grid. The total mass can be
computed by the following Riemann type sum
\begin{equation*}
\mu _{n}=\Delta x\sum_{j=1}^{J-1}U_{j}^{n}
\end{equation*}
A discrete version of (\ref{5.2}) may look like
\begin{eqnarray}
\frac{\mu _{n+1}-\mu _{n}}{\Delta t} &=&\Delta x\sum_{j=1}^{J-1}\frac{%
U_{j}^{n+1}-U_{j}^{n}}{\Delta t},  \label{5.3} \\
&=&\Delta x\sum_{j=1}^{J-1}\frac{U_{j-1}^{n+1}-2U_{j}^{n+1}+U_{j+1}^{n+1}}{%
(\Delta x)^{2}},  \notag \\
&=&\frac{U_{J}^{n+1}-U_{J-1}^{n+1}}{\Delta x}-\frac{U_{1}^{n+1}-U_{0}^{n+1}}{%
\Delta x},  \notag \\
&=&2\phi (\tau _{n+1}).  \notag
\end{eqnarray}
Next we shall explicitly find $\mu _{n}$ and the time switches $t_{n}$ with
some restriction on the mesh size. In the first subinterval $0<\tau
_{n}<t_{1},$ we choose $\Delta t=\Delta t_{0}=\dfrac{M}{2N_{0}},$ where $%
N_{0}$ is a positive integer, along with the flux $\phi (\tau _{n})=1.$ We
obtain
\begin{equation*}
\mu _{n}=2n\Delta t_{0}=2\tau _{n},\,\,\,\,\,\,0\le \tau _{n}\le t_{1}
\end{equation*}
This special choice for $\Delta t_{0}$ leads to
\begin{equation*}
\mu _{N_{0}}=M
\end{equation*}
i.e. $N_{0}$ specifies the first time switch
\begin{equation*}
t_{1}=\tau _{N_{0}}=\frac{M}{2}
\end{equation*}
which coincides with the one computed analytically. Next, we reverse the
flux, i.e. $\phi (\tau _{n})=-1,$ for $n>N_{0},$ and define $\Delta t=\dfrac{%
M-m}{2N},$ for some positive integer $N.$ Then in view of
(\ref{5.3}), we obtain
\begin{equation*}
\mu _{n}=M-2(n-N_{0})\Delta t,\,\,\,\,\,\,\,n\ge N_{0}.
\end{equation*}
Notice that in the second stage, the time mesh size is not
necessarily equal to the one at the first stage. With a carefully
chosen $\Delta t,$ we can easily get
\begin{equation*}
\mu _{N_{0}+N}=m
\end{equation*}
with the second time switch equal to
\begin{eqnarray*}
t_{2} &=&N_{0}\Delta t_{0}+N\Delta t \\
&=&t_{1}+(\frac{M}{2}-\frac{N}{2}).
\end{eqnarray*}
This result is in agreement with the one computed analytically.

For the next upcoming stages, we keep the time mesh size as $\Delta t=\dfrac{%
M-m}{2N}.$ The following are the total mass and the time switches for the
third and the forth stages respectively
\begin{eqnarray*}
\mu _{n} &=&m+2(n-N_{0}-N)\Delta t,\,\,\,\,\,\,\,n\ge N_{0}+N, \\
&& \\
t_{3} &=&N_{0}\Delta t_{0}+2N\Delta t, \\
&=&t_{2}+(\frac{M}{2}-\frac{N}{2}), \\
&& \\
\mu _{n} &=&M-2(n-N_{0}-2N)\Delta t,\,\,\,\,\,\,\,n\ge N_{0}+2N, \\
&& \\
t_{4} &=&N_{0}\Delta t_{0}+3N\Delta t, \\
&=&t_{3}+(\frac{M}{2}-\frac{N}{2}),
\end{eqnarray*}
This inductively gives
\begin{eqnarray*}
\mu _{n} &=&\left\{
\begin{array}{lll}
M-2(n-N_{0}-kN)\Delta t, & \text{if} & N_{0}+kN\le n\le
N_{0}+(k+1)N,\,\,\,\,k\text{ is even} \\
&  &  \\
m+2(n-N_{0}-kN)\Delta t, & \text{if} & N_{0}+kN\le n\le
N_{0}+(k+1)N,\,\,\,\,k\text{ is odd}
\end{array}
\right.  \\
&& \\
t_{k} &=&N_{0}\Delta t_{0}+(k-1)N\Delta t, \\
&=&t_{k-1}+(\frac{M}{2}-\frac{N}{2}).
\end{eqnarray*}
The above calculations shows that if the time grid is chosen properly, the
time switches computed by the above difference scheme coincides with those
computed analytically.

\noindent \textbf{Remark:} The situation when the $\tau _{n}$ is
equispaced does not, in general, generate the exact time switches
$t_{n}$. However, if we take the switching criteria as, $\mu _{n}\ge
M,$ and $\mu _{n}\le m$ instead of $\mu _{n}=M$ and $\mu _{n}=m,$
respectively, where $n$ is the least integer that satisfies such
inequalities, we obtain a new set of switching points, say $T_{n}$,
with the errors
\begin{eqnarray*}
0 &\le &T_{1}-t_{1}<\Delta t, \\
0 &\le &T_{2}-t_{2}<2\Delta t, \\
0 &\le &T_{3}-t_{3}<3\Delta t, \\
&&\vdots  \\
0 &\le &T_{k}-t_{k}<k\Delta t.
\end{eqnarray*}
If the maximum time limit $T$ is specified, with $\Delta t=T/N,$ for
some integer $N$, then for a fixed integer $k,$ we will have
$T_{k}-t_{k}<k\Delta t=kT/N$ convergent to $0$ as $N\rightarrow
\infty .$

\section{Numerical Example}
In this section, we consider a finite difference method to
discretize the problem

\begin{eqnarray}
u_t=\al u_{xx}, \quad 0<x<1, \quad 0<t\le T\notag\\
-u_x(0,t)=u_x(1,t)=\phi(t), \quad 0<t\le T\notag\\
u(x,0)=0, \quad 0<x<1 \notag
\end{eqnarray}
where the boundary control function is \bee
\phi(t)= \begin{cases} 1, \quad &t_{2n} \leq t\leq t_{2n+1}\\
-1, & \text{elsewhere} \end{cases} \notag \ene and $\{ t_n\}$
depends on \bee \mu (t)=\int_0^1 u(x,t)dx \notag \ene where
\begin{alignat*}{2}
\mu (t_{2n} )&= 0.1,  &\quad & n=1,2, \ldots\\
\mu (t_{2n+1} ) &=0.2,  & & n=0,1, \ldots.
\end{alignat*}
The time limit and the diffusivity constant are taken as $T=10$
and $\al =0.05$.

Let's consider the space and time discretization

i) $\Delta x=\frac{1}{J} , \quad x_j=j\Delta x, \quad j=0,1,\dots
,J$

ii) $\Delta t=\frac{T}{N} , \quad \tau_n =n\Delta t, \quad
n=0,\dots ,N$

\noindent where $J=50$ and $N=200$ . The integer $N$ is chosen
large enough so that the time step $\Delta t$ is much smaller than
an estimated differences between two consecutive values of the
time switches.

We consider the backward implicit finite difference scheme
$$
\frac{U_j^{n+1}-U_j^{n}}{\Delta t} = \al \frac{U_{j-1}^{n+1}
-2U_j^{n+1} +U_{j+1}^{n+1}}{(\Delta x)^2}
$$
which can be written as \bee -\nu U_{j-1}^{n+1} +(1+2\nu )U_j^{n+1}
-\nu U_{j+1}^{n+1} =U_j^n  \ene where $\nu =\al\Delta t/(\Delta
x)^2$, $j=1,\dots ,J-1$ and $n=0,1,\dots, N-1$. The initial
condition is $U_j^0=0$ for $j=1,\dots ,J-1$, and the boundary
conditions are \bee -\frac{U_1^n-U_0^n}{\Delta x}
=\frac{U_J^n-U_{J-1}^n}{\Delta x} =\phi(\tau_n )\ene for
$n=0,1,\dots ,N$. The total mass integral is calculated by the
following trapezoidal rule \bee \mu_n =\frac{h}{2} \sum_{j=0}^{N-1}
\left( U_j^{n+1} +U_{j+1}^n \right). \label{5.11}\ene The numerical
experiment is carried out in the following way. We start by setting
the flux at $\phi=1$ then we solve a tridiagonal system coming out
of the difference method. We evaluate the total mass $\mu_n$ and
compare it with the upper threshold $M=0.2$. We move to the next
time step while keeping the flux at $\phi =1$, as long as $\mu_n<M$,
otherwise, we switch the flux to $\phi=-1$, as long as $\mu_n\ge M $
. At the moment, say $\tau_{n_1}$, for some integer $n_1$, when the
total mass exceeds $M$ for the first time, we take $T_1=\tau_{n_1}$
as an approximation for the first time switch. With $\phi=0$, we
proceed our solution along the time, as long as $\mu_n$ does not
fall below the threshold $m=0.1$. By the moment, when $\mu_{n_2} \le
m$, for some integer $n_2$, we set $T_2=\tau_{n_2}$, and we switch
the flux back to $\phi =1$ at the next step. We keep switching the
flux between ($\phi=1$) and ($\phi=-1$) and calculating the time
switches $T_k$ until the end of the run when $\tau_n = 10$.

Table (1) shows the times switches $T_n$. As we can see there, the
difference between any two consecutive time switches tend to 0.95.
For the same set of data, graphs (1) through (4) show the
concentration versus the space at consecutive time steps. The
graphs are obtained for different stages, where at each stage the
flux is kept constant at the end points.  A profile of the
concentrations at $x=0.5$ for various times is shown  in graph (5)
with the same specified data. Graph (6) shows the total mass
computed analytically and numerically .

\begin{table}
\begin{tabular}{|c|c|c|}
\hline
  $n$ &  $T_n$ & $T_n - T_{n-1}$ \\\hline
1&1.9500&1.9000    \\
2&2.9000&0.9500    \\
3&3.8500&0.9500    \\
4&4.8000&0.9500    \\
5&5.7500&0.9500    \\
6&6.7000&0.9500    \\
7&7.6500&0.9500    \\
8&8.6000&0.9500    \\
9&9.5500&0.9500    \\

 \hline
\end{tabular}
\medskip \caption{The Time switches $T_n$ and the differences $T_n-T_{n-1}$.
 Note the differences between any two consecutive times tends to 0.95.}
\end{table}

\begin{figure}[p]
\scalebox{0.8}{\includegraphics{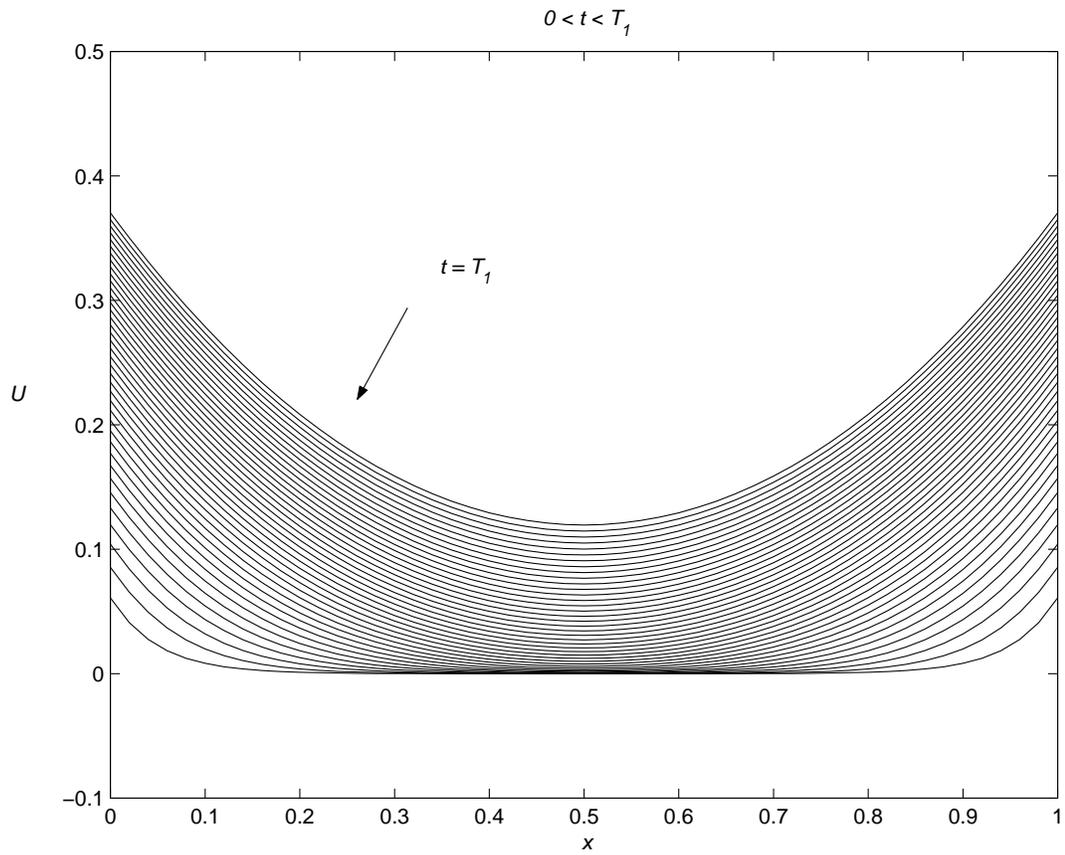}}
 \caption{The first
stage where the flux $\phi$ is held at 1 at the end points. Each
curve shows the concentration profile at various discrete time
steps $\tau_n=n \Delta t$. As the time goes on, the level of
concentrations gets higher }
\end{figure}


\begin{figure}[p]
\scalebox{0.8}{\includegraphics{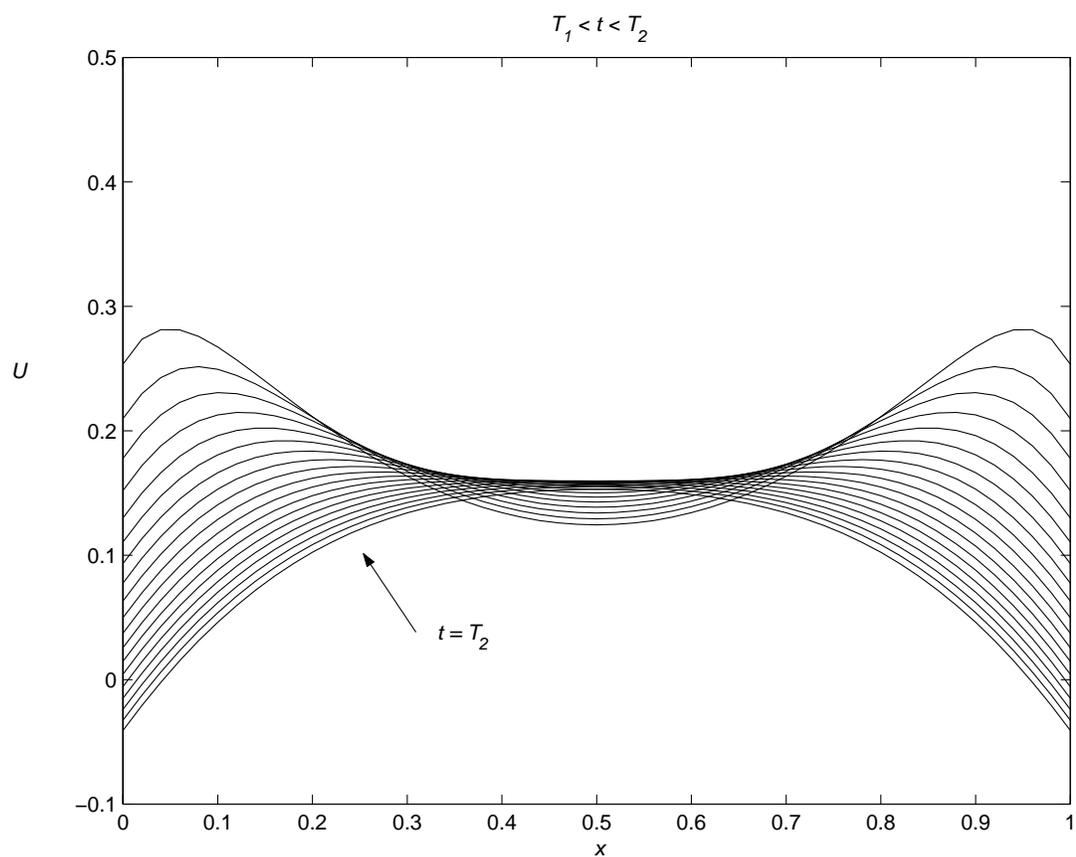}} \caption{The second
stage where the flux $\phi$ is held at -1 at the end points.
 As the time goes on, the level of
concentrations decreases. Notice the fluctuations when the
concentration is dropped suddenly to 0 at the beginning of the
stage.}
\end{figure}

\begin{figure}[p]
\scalebox{0.8}{\includegraphics{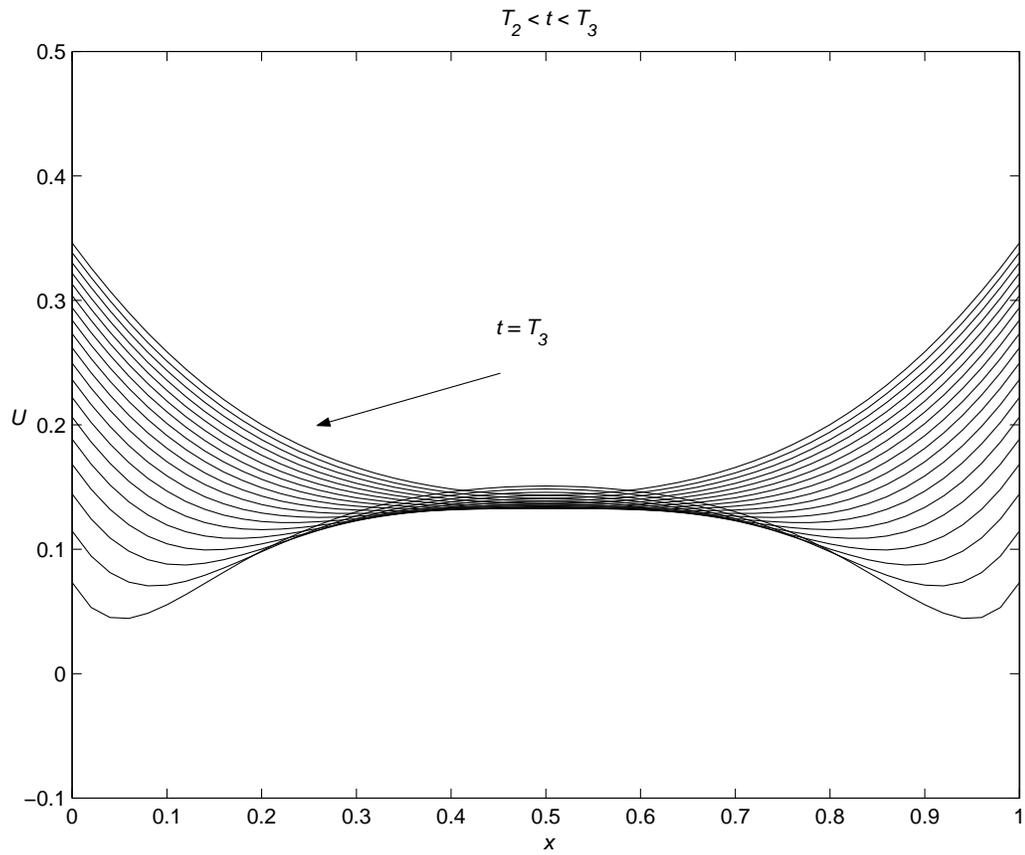}}

\caption{The third stage where the flux $\phi$ is switched to 1 at
the end points. Each curve shows the concentration profile at
various discrete time steps. Notice the fluctuations due to the
sudden change on the concentrations. After a little while, the
concentrations levels increase monotonically. }
\end{figure}

\begin{figure}[p]
\scalebox{0.8}{\includegraphics{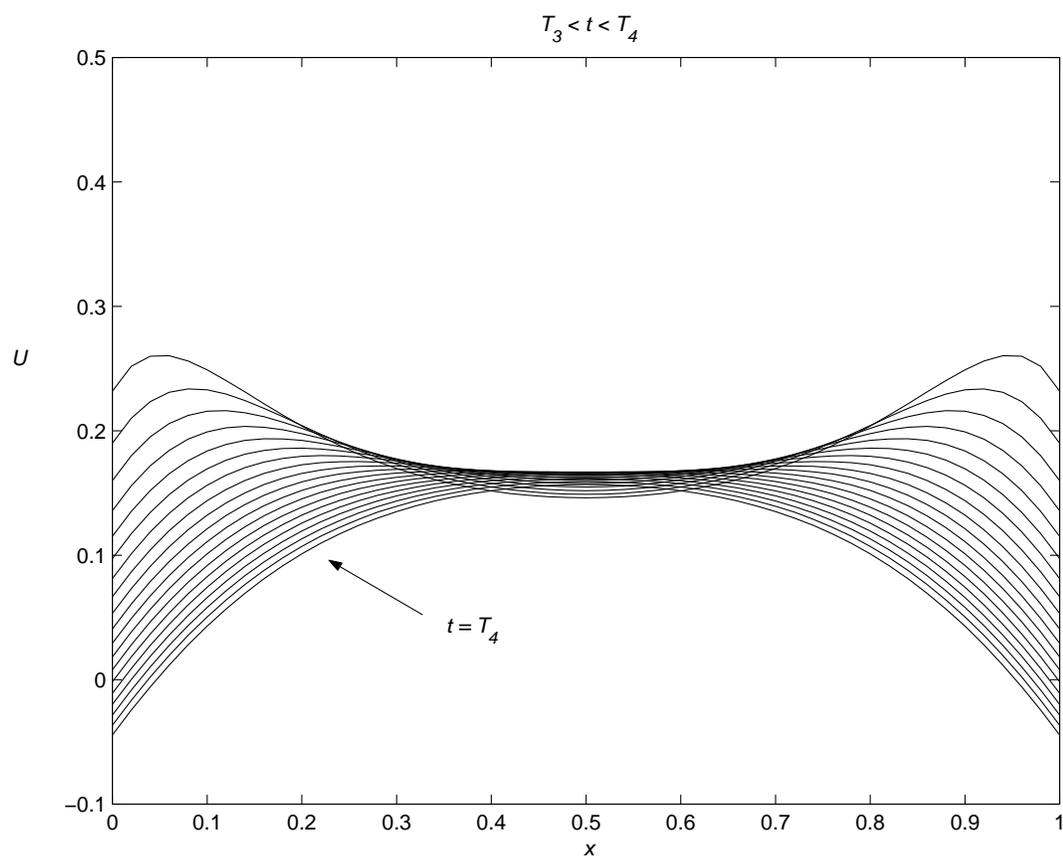}}

\caption{The fourth stage where the flux $\phi$ is switched to -1
at the end points. Notice the similarity with the second stage.}
\end{figure}

\begin{figure}[p]
\scalebox{0.8}{\includegraphics{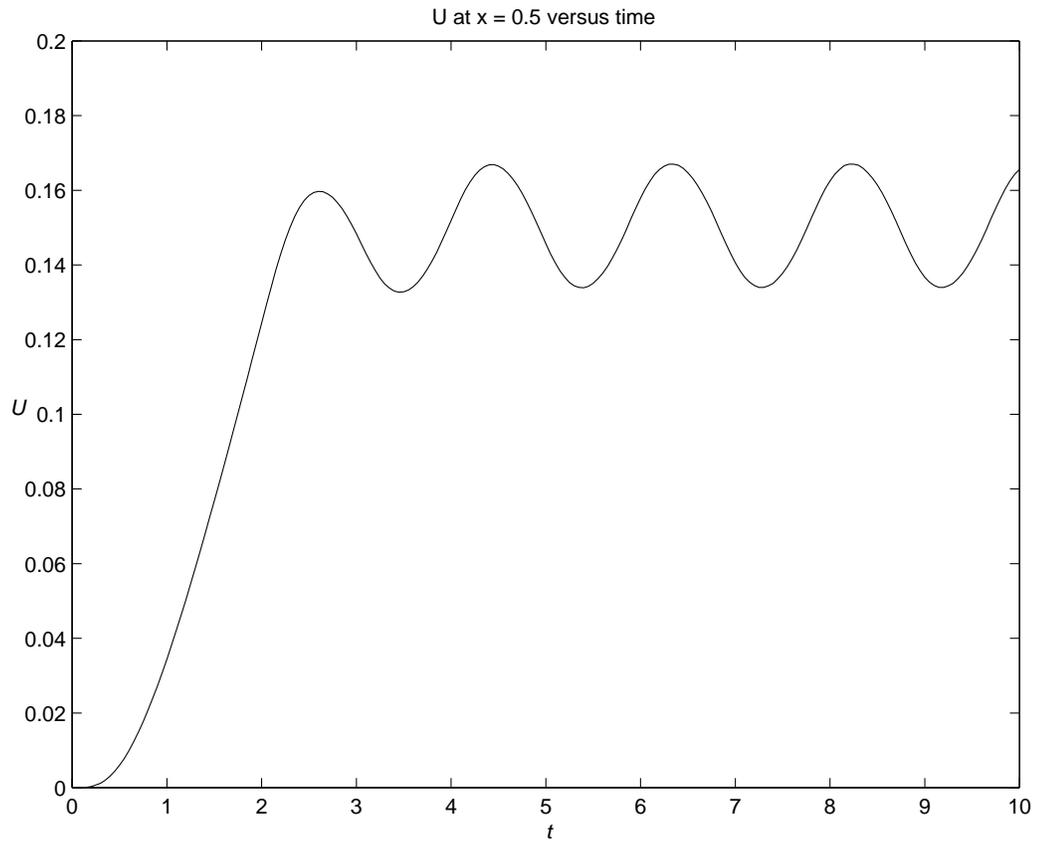}}

\caption{ The concentration profile $U$ at $x=0.5$ versus the time
shows periodic behavior due to the periodic change of the boundary
conditions}
\end{figure}

\begin{figure}[p]
\scalebox{0.8}{\includegraphics{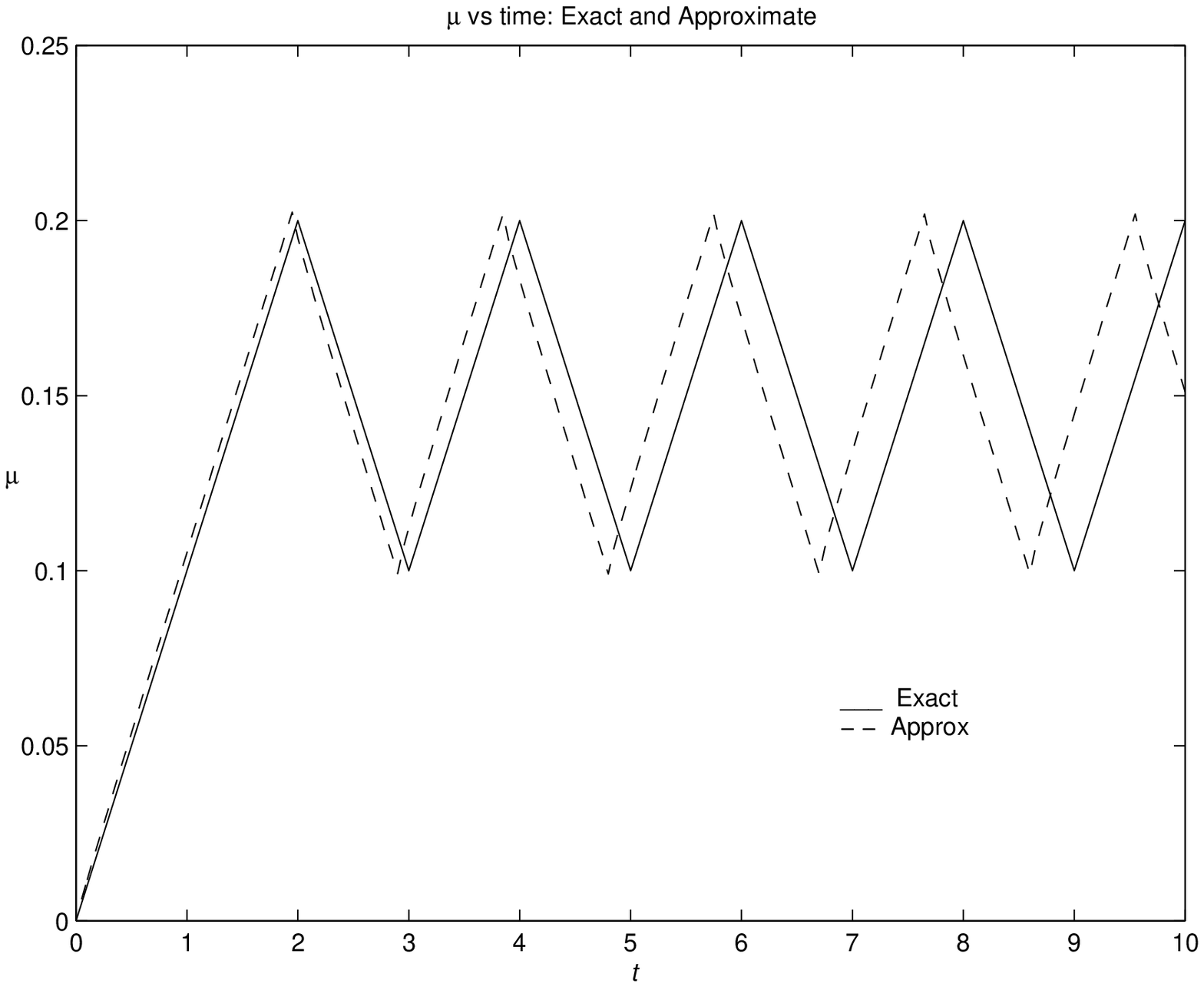}}

\caption{The total mass computed analytically and numerically.
Note how the error in calculating the time switches accumulates as
the time gets large.}
\end{figure}


\begin{thebibliography}{99}
\bibitem{can1} J. R. Cannon, The one-dimensional heat equation,
{\it Encyclopedia of Mathematics and its Applications} {\bf 23}
(1984).

\bibitem{can2} J. R. Cannon, Yapping Lin and Shuzhan Xu, Numerical
procedures for the determination of an unknown coefficient in
semi-linear parabolic differential equation. {\it Inverse
Problems} {\bf 10} (1994), 227-243.

\bibitem{canv1} J.R. Cannon and J. Van der Hoek, The existence and
a continuous dependance result for the solution of the heat
equation subject to the specification of energy. {\it Boll. Uni.
Math. Ital. Suppl.} {bf 1}(1981), 253-282.

\bibitem{canv2} J.R. Cannon and J. Van der Hoek, An emplicit
finite difference scheme for the diffusion of mass in a portion of
a domain. {\it Numerical Solutions of Partial Differential
Equations.} (J. Noye, ed) ,527-539, {\it North-Holland,
Amsterdam,} 1982.

\bibitem{canv3} J.R. Cannon and J. Van der Hoek, Diffusion subject
to a specification of mass. {\it J. Math. Anal. Appl.} {\bf 115}
(1986), No.2, 517-529.

\bibitem{can05} J.R.\ Cannon \& M.\ Salman, The Utilization of Total Mass to
Determine the Switching Points in the Symmetric Boundary Control
Problem with a Linear Reaction Term. \textit{Journal of Mathematical
Analysis and Applications} {\bf 311} (2005), No.\ 1, 147--161.

\bibitem{09} J.R.\ Cannon \& M.\ Salman,  A Boundary Control Problem with
a Nonlinear Reaction Term.  \textit{Electron. J. Diff. Eqns.}, Conf
17 (2009), pp. 39-49.

\bibitem{day1} W. A. Day, Existence of a property of solutions of
the heat equation to linear thermoelasticity and other theories.
{\it Quart. Appl. Math.}, {\bf 40}, (1982) 319-330.


\bibitem{day2} W. A. Day, A decreasing property of solutions of
a parabolic equation with applications to thermoelasticity. {\it
Quart. Appl. Math.}, {\bf 41}, (1983) 468-4475.




\bibitem{kev} J. Kevorkian, Partial differential equations,
analytical solution techniques. 2nd edition. {\it Texts in Applied Mathematics} {\bf 35} (2000).

\bibitem{pao1} Pao, C. V. Blowing-up of solution for a nonlocal reaction-diffusion
problem in combustion theory. {\it J. Math. Anal. Appl.} {\bf 166} (1992), no. 2, 591-600.

\bibitem{pao2} Pao, C. V. Dynamics of reaction-diffusion equations
with nonlocal boundary conditions. {\it Quart. Appl. Math.} {\bf 53} (1995), no. 1, 173-186.

\bibitem{pao3} Pao, C. V. Reaction diffusion equations with nonlocal
boundary and nonlocal initial conditions. {\it J. Math. Anal. Appl.} {\bf 195} (1995), no. 3, 702-718.

\bibitem{pao4} Pao, C. V. Dynamics of weakly coupled parabolic systems with
nonlocal boundary conditions. {\it Advances in nonlinear dynamics,} 319-327,
Stability Control Theory Methods Appl., 5, {\it Gordon and Breach, Amsterdam,} 1997.

\bibitem{pao5} Pao, C. V. Asymptotic behavior of solutions of reaction-diffusion
equations with nonlocal boundary conditions. Positive solutions of nonlinear problems.
{\it J. Comput. Appl. Math.} {\bf 88} (1998), no. 1, 225-238.

\bibitem{pao6} Pao, C. V. Numerical solutions of reaction-diffusion equations
 with nonlocal boundary conditions. {\it J. Comput. Appl. Math.} {\bf 136} (2001), no. 1-2, 227-243.

\bibitem{prot} Protter, Murray H.; Weinberger, Hans F. Maximum principles in differential equations.
{\it Prentice-Hall, Inc., Englewood Cliffs, N.J.} 1967 x+261 pp.

\bibitem{can}  M.\ Salman, The utilization of total mass
to determine the switching points in the symmetric boundary control
of a diffusion problem. (Manuscript)

\end{thebibliography}
\end{document}